# Statistical models, likelihood, penalized likelihood and hierarchical likelihood

Daniel Commenges

Epidemiology and Biostatistics Research Center, INSERM

Université Victor Segalen Bordeaux 2

146 rue Léo Saignat, Bordeaux, 33076, France

Tel: (33) 5 57 57 11 82; Fax (33) 5 56 24 00 81

November 2, 2018

We give an overview of statistical models and likelihood, together with two of its variants: penalized and hierarchical likelihood. The Kullback-Leibler divergence is referred to repeatedly, for defining the misspecification risk of a model, for grounding the likelihood and the likelihood crossvalidation which can be used for choosing weights in penalized likelihood. Families of penalized likelihood and sieves estimators are shown to be equivalent. The similarity of these likelihood with a posteriori distributions in a Bayesian approach is considered.

# 1 Introduction

Since its proposal by Fisher (1922), likelihood inference has occupied a central position in statistical inference. In some situations modified versions of the likelihood have been proposed. Marginal, conditional, profile and partial likelihoods have been proposed to get rid of nuisance parameters. Pseudo-likelihood and hierarchical likelihood may be used to circumvent numerical problems in the computation of the likelihood, generally due to multiple integrals. Penalized likelihood has been proposed to introduce a smoothness a priori knowledge on functions, thus leading to smooth estimators. Several review have already been proposed, for instance Lee and Nelder (1992), but it is nearly impossible in a single paper to describe with some details all the types of likelihoods that have been proposed. This paper aims at describing the conventional likelihood and two of its variants: penalized and hierarchical likelihoods. The aim of this paper is not to give the properties of the estimators obtained by maximizing these likelihood but rather to describe these three likelihoods together with their link to the Kullback-Leibler divergence. This interest more turned to the foundations than to the properties, leads us to first develop some reflexions and definitions about statistical models and to give a slightly extended version of the Kullback-Leibler divergence.

In section 2 we recall the definition of a density and the relationship be-tween a density in the sample space and for a random variable. In section 3 we give a slightly extended version of the Kullbaclk-Leibler divergence (making it explicit that it also depends on a sigma-field). Section 4 gives an account of statistical models, distinguishing mere statistical families from statistical models, and defining the misspecification risk. Section 5 presents the likelihood and discusses issues about its computation and the performance of the estimator of the maximum likelihood in terms of Kullback-Leibler risk. In section 6 we define the penalized likelihood and show that for a family of penalized likelihood estimators there is an identical family of sieves estimators. In section 7 we describe the hierarchical likelihood. In section 8 we briefly sketch the possible unification of these likelihoods through a Bayesian representation allowing to consider the maximum (possibly penalized) likelihood estimators as MAP estimators; this question however cannot be easily settled due to the non-invariance of the MAP for reparametrization. There is a short conclusion.

## 2 Definition of a density

Consider a measurable space $(\mathcal{S}, \mathbf{A})$ and two measures $\mu$ and $\nu$ with $\mu$ absolutely continuous relatively to $\nu$. For $\mathcal{G}$ a sub-$\sigma$-field of $\mathcal{F}$ the Radon-Nikodym derivative of $\mu$ with respect to $\nu$ on $\mathcal{X}$, denoted by: $\frac{d\mu}{d\nu}_{|\mathcal{G}}$ is the $\mathcal{G}$-measurable random variable such that

$$\mu(F) = \int_G \frac{d\mu}{d\nu}_{|\mathcal{G}} d\nu, G \in \mathcal{G}.$$

The Radon-Nikodym derivative is also called the density. We are interested in the case where $\mu$ is a probability measure, that we will call $P^1$; $\nu$ may also be a probability measure, $P^0$. In that case we can speak of likelihood ratio and denote it $\mathcal{L}_{\mathcal{G}}^{P^1/P^0}$. In order to speak of likelihood function, we have to define a model (see section 4). Note that likelihood ratios (as Radon-Nykodym derivatives) are defined with respect to a sigma-field. If $\mathcal{H}$ and $\mathcal{G}$ are different sigma-fields, $\frac{dP^1}{dP^0}_{|\mathcal{H}}$ and $\frac{dP^1}{dP^0}_{|\mathcal{G}}$ are different, but if $\mathcal{H} \subset \mathcal{G}$ the former can be expressed as a conditional expectation (given $\mathcal{H}$) of the latter and we have the fundamental formula:

$$\frac{dP^1}{dP^0}_{|\mathcal{H}} = \mathrm{E}_{P^0}\left[\frac{dP^1}{dP^0}_{|\mathcal{G}}|\mathcal{H}\right].$$

Consider now the case where the measurable space $(\Omega, \mathcal{F})$ is the sample space of an experiment and define a random variable $X$, a measurable function from $(\Omega, \mathcal{F})$ to $(\Re, \mathcal{B})$. We shall write in bold character a probability on $(\Omega, \mathcal{F})$, for instance, $\boldsymbol{P}^1$. The couple $(\boldsymbol{P}^1, X)$ induces a probability measure on $(\Re, \mathcal{B})$ defined by: $P_X^1(B) = \boldsymbol{P}^1 o X^{-1}(B), B \in \mathcal{B}$. This probability measure is called the distribution of $X$. If this probability measure is absolutely continuous with respect to Lebesgue (resp. counting) measure, one speaks of continuous (resp. discrete) variable. For instance, for a continuous variable we define the density $f_X^1 = \frac{dP_X^1}{d\lambda}$, which is the usual probability density function (p.d.f.). Note that the p.d.f. depends on both $\boldsymbol{P}^1$ and $X$, while $\frac{d\boldsymbol{P}^1}{d\boldsymbol{P}^0}_{|\mathcal{X}}$ depends on $\mathcal{X}$ but not on a specific random variable $X$. Often in applied statistics one works only with distributions, but this may let some problems unsolved.

**Example 1**. Consider the case where concentrations of CD4 lymphocytes are measured. $\Omega$ represents the set of physical concentrations that may happen. Let the random variables $X$ and $Y$ express the concentration in number of CD4 by mm³ and by $ml$ respectively. Thus we have $Y = 10^3 X$.

So $X$ and $Y$ are different, although they are informationally equivalent. For instance the events $\{\omega : X(\omega) = 400\}$ and $\{\omega : Y(\omega) = 400000\}$ are the same. The densities of $X$ and $Y$, for the same $\boldsymbol{P}^1$ on $(\Omega, \mathcal{F})$, are obviously different. So, if we look only at distributions, we shall have difficulties to define rigorously what a model is.

## 3 The Kullback-Leibler risk

Many problems in statistical inference can be treated from the point of view of decision theory. That is, estimators for instance are chosen as minimizing some risk function. The most important risk function is based on the Kullback-Leibler divergence: maximum likelihood estimators, use of Akaike criterion or likelihood crossvalidation can be gounded on the Kullback-Leibler divergence. Given a probability $\boldsymbol{P}^2$ absolutely continuous with respect to a probability $\boldsymbol{P}^1$ and $\mathcal{X}$ a sub-$\sigma$-field of $\mathcal{F}$, the loss using $\boldsymbol{P}^2$ in place of $\boldsymbol{P}^1$ is the log-likelihood ratio $L_{\mathcal{X}}^{\boldsymbol{P}^1/\boldsymbol{P}^2} = \log \frac{d\boldsymbol{P}^1}{d\boldsymbol{P}^2}_{|\mathcal{X}}$. Its expectation is $\mathrm{E}_{\boldsymbol{P}^1}[L_{\mathcal{X}}^{\boldsymbol{P}^1/\boldsymbol{P}^2}]$. This is the Kullback-Leibler risk, also called divergence (Kullback and Leibler, 1951; Kullback, 1959) or information deviation (Cencov, 1972) or entropy (Akaike, 1973). The different names of this quantity reflects its central position in statistical theory, being connected to several fields of the theory. Several notations have been used by different authors. Here we choose the Cencov notation:

$$\boldsymbol{I}(\boldsymbol{P}^2|\boldsymbol{P}^1; \mathcal{X}) = \mathrm{E}_{\boldsymbol{P}^1}[L_{\mathcal{X}}^{\boldsymbol{P}^1/\boldsymbol{P}^2}]$$

If $\mathcal{X}$ is the largest sigma-field defined on the space we omit it in the notation. Note that the Kullback-Leibler risk is asymmetric and hence does not define a distance between probabilities; we have to take on this fact. If $X$ is a random variable with p.d.f. $f_X^1$ and $f_X^2$ under $\boldsymbol{P}^1$ and $\boldsymbol{P}^2$ respectively we have $\frac{d\boldsymbol{P}^1}{d\boldsymbol{P}^2}_{|\mathcal{X}} = \frac{f_X^1(X)}{f_X^2(X)}$ and the divergence of the distribution $P_X^2$ relative to $P_X^1$ can be written:

$$\mathcal{I}(P_X^2|P_X^1) = \int \log \frac{f_X^1(x)}{f_X^2(x)} f_X^1(x) dx. \tag{1}$$

We have that $\boldsymbol{I}(\boldsymbol{P}^2|\boldsymbol{P}^1; \mathcal{X}) = \mathcal{I}(P_X^2|P_X^1)$. Note that on $(\Omega, \mathcal{F})$ we have to specify that we assess the divergence on $\mathcal{X}$; we might assess it on a different sigma-field and would of course obtain a different result. This gives more flexibility. In particular we shall use this in the case of incomplete data. The observation is represented by a sigma-field $\mathcal{O}$. Suppose we are interested to make inference about the true probability on $\mathcal{X}$. We have complete data if our observation is $\mathcal{O} = \mathcal{X}$. With incomplete data, in the case where the mechanism leading to incomplete data is deterministic, we have $\mathcal{O} \subset \mathcal{X}$. In that case it will be very difficult to estimate $\boldsymbol{I}(\boldsymbol{P}^2|\boldsymbol{P}^1; \mathcal{X})$ and it will be more realistic to use $\boldsymbol{I}(\boldsymbol{P}^2|\boldsymbol{P}^1; \mathcal{O}) = \mathrm{E}_{\boldsymbol{P}^1}[L_\mathcal{O}^{\boldsymbol{P}^1/\boldsymbol{P}^2}]$. We need this flexibility to extend Akaike's argument for the likelihood and for developing model choice criteria to situations with incomplete data.

**Example 2**. Suppose we are interested in modeling the time to an event, $X$, and we wish to evaluate the divergence of $\boldsymbol{P}^2$ with respect to $\boldsymbol{P}^1$. It is natural to compute the divergence on the sigma-field $\mathcal{X}$ generated by $X$, $\boldsymbol{I}(\boldsymbol{P}^2|\boldsymbol{P}^1; \mathcal{X}) = \mathcal{I}(P_X^2|P_X^1)$ given by formula (1). Suppose however that observation of $X$ is right-censored at a fixed time $C$. We observe $(\tilde{X}, \delta)$ where $\tilde{X} = \min(X, C)$ and $\delta = 1_{\{X \leq C\}}$. Thus on $\{X \leq C\}$ we observe all the events of $\mathcal{X}$ but on $\{X > C\}$ we observe no more detailed event. If we represent the observation by the sigma-field $\mathcal{O}$ we can say that $\mathcal{O}$ is generated by $(\tilde{X}, \delta)$. It is clear that we have $\mathcal{O} \subset \mathcal{X}$. Although in theory it is still interesting to compute the divergence of $\boldsymbol{P}^2$ with respect to $\boldsymbol{P}^1$ on the sigma-field $\mathcal{X}$ it is also interesting to compute it on the observed sigma-field, that is $\boldsymbol{I}(\boldsymbol{P}^2|\boldsymbol{P}^1;\mathcal{O})$.

It can be proved by simple probabilistic arguments that on $\{X \leq C\}$ we have $\frac{d\boldsymbol{P}^1}{d\boldsymbol{P}^2}\big|_{\mathcal{O}} = \frac{f_X^1(X)}{f_X^2(X)}$ and on $\{X > C\}$ we have $\frac{d\boldsymbol{P}^1}{d\boldsymbol{P}^2}\big|_{\mathcal{O}} = \frac{S_X^1(C)}{S_X^2(C)}$ and thus

$$\boldsymbol{I}(\boldsymbol{P}^2|\boldsymbol{P}^1;\mathcal{O}) = \int_0^C \log \frac{f_X^1(x)}{f_X^2(x)} \, f_X^1(x) dx + \log \frac{S_X^1(C)}{S_X^2(C)} \, S^1(C).$$

This will take all its importance in section 5 where $\boldsymbol{P}^1$ will be the true unknown probability (denoted $\boldsymbol{P}^*$); the problem will not be to compute but to estimate this divergence.

## 4 Statistical models and families

### 4.1 Statistical families

Consider a measure space $(\mathcal{S}, \boldsymbol{A}, \mu)$. We consider a subset $\mathcal{P}$ of the probabilities on $(\mathcal{S}, \boldsymbol{A}, \mu)$. We shall call such a subset a family of probabilities. We may parametrize this family. Following Hoffmann-Jorgensen (1994) a parametrization can be represented by a function from a set $\Theta$ with values in $\mathcal{P}$: $\theta \to P^\theta$. It is desirable that this function be one-to-one, a property linked to the identifiability issue which will be discussed later in this section. The parametrization associated to the family of probabilities $\mathcal{P}$ can be denoted $\Pi = (P^\theta; \theta \in \Theta)$ and we have $\mathcal{P} = \{P^\theta; \theta \in \Theta\}$. We may denote $\Pi \sim \mathcal{P}$. If $\Pi_1 \sim \mathcal{P}$ and $\Pi_2 \sim \mathcal{P}$, $\Pi_1$ and $\Pi_2$ are two parameterizations of the same family of probabilities and we may note $\Pi_1 \sim \Pi_2$.

However we do not consider that a parametrized family on $(\Re, \mathcal{B})$ representing a family of distributions of a random variable is sufficient to specify a statistical model (here, we do not follow Hoffmann-Jorgensen, 1994). This is because the distributions depend on the random variables chosen, as exemplified in section 2.

### 4.2 Statistical models

A family of probabilities on the sample space of an experiment $(\Omega, \mathcal{F})$ will be called a statistical model and a parametrization of this family will be called a parametrized statistical model.

**Definition 1** *Two parametrized statistical models* $\mathbf{\Pi} = (\boldsymbol{P}^\theta, \theta \in \Theta)$ *on* $\mathcal{X}$ *and* $\mathbf{\Pi}' = (\boldsymbol{P}^\gamma, \gamma \in \Gamma)$ *on* $\mathcal{Y}$ *are equivalent (in the sense that they specify the same statistical model) if* $\mathcal{X} = \mathcal{Y}$ *and they specify the same family of probability on* $(\Omega, \mathcal{X})$.

The couple $(X, \mathbf{\Pi})$ of a random variable and a parametrized statistical model induces the parametrized family (of distributions) on $(\Re, \mathcal{B})$: $\Pi_X = (P_X^\theta; \theta \in \Theta)$. Conversely the couple $(X, \Pi_X)$ induces $\mathbf{\Pi}$ if $\mathcal{X} = \mathcal{F}$. In that case we may describe the statistical model by $(X, \Pi_X)$. Two different random variables $X$ and $Y$ induce two (generally different) parametrized families on $(\Re, \mathcal{B})$, $\Pi_X$ and $\Pi_Y$. Conversely one may ask whether the couples $(X, \Pi_X)$ and $(Y, \Pi_Y)$ define equal or equivalent parametrized statistical models. We need the definition of "informationally equivalent" random variables (or more generally random elements).

**Definition 2** *X and Y are informationally equivalent if the sigma-fields $\mathcal{X}$ and $\mathcal{Y}$ generated by X and Y are equal.*

Each couple $(X, P_X^\theta)$ induces a probability on $(\Omega, \mathcal{X})$ $\boldsymbol{P}^{X,\theta} = P_X^\theta oX$ and thus the couple $(X, \Pi_X)$ induces the parametrized statistical model $(\boldsymbol{P}^{X,\theta}, \theta \in \Theta)$. Similarly each couple $(Y, P_Y^\gamma)$ induces a probability on $(\Omega, \mathcal{Y})$ $\boldsymbol{P}^{Y,\gamma} = P_Y^\gamma oY$ and the couple $(Y, \Pi_Y)$ induces the parametrized statistical model $(\boldsymbol{P}^{Y,\gamma}, \gamma \in \Gamma)$. Tautologically we will say that $(X, \Pi_X)$ and $(Y, \Pi_Y)$ define the same statistical models if $(\boldsymbol{P}^{X,\theta}, \theta \in \Theta)$ and $(\boldsymbol{P}^{Y,\gamma}, \gamma \in \Gamma)$ are equivalent.

**Example 1 continued**

(i) $\Pi_X = (\mathcal{N}(10^3; \sigma^2), \sigma^2 > 0)$ and $\Pi_Y = (\mathcal{N}(10^3; \sigma^2), \sigma^2 > 0)$ are the same parametrized families on $(\Re, \mathcal{B})$. However if X and Y are measurements of the same quantity in different units, these parametrized families correspond to different statistical models.

(ii) $\Pi_X = (\mathcal{N}(\mu, \sigma^2); \mu \in \Re, \sigma^2 > 0)$ and $\Pi_Y = (\mathcal{N}(\mu, \sigma^2); \mu \in \Re, \sigma^2 > 0)$ are the same parametrized family on $(\Re, \mathcal{B})$. $(X, \Pi_X)$ and $(Y, \Pi_Y)$ specify the same statistical model but not the same parametrized statistical model.

(iii) $\Pi_X = (\mathcal{N}(10^3; \sigma^2), \sigma^2 > 0)$ and $\Pi_Y = (\mathcal{N}(10^6; 10^6\sigma^2), \sigma^2 > 0)$ are different families on $(\Re, \mathcal{B})$. However $(X, \Pi_X)$ and $(Y, \Pi_Y)$ specify the same statistical model (with the same parametrization).

For sake of simplicity we have considered distributions of real random variables. The same can be said about random variables with values in $\Re^d$ or stochastic processes which are random elements with values in a Skorohod space. Commenges and Gégout-Petit (2007) gave an instance of two informationally equivalent processes. The events described by an irreversible three-state process $X = (X_t)$, where $X_t$ takes values $0, 1, 2$, can be described by a bivariate counting process $N = (N_1, N_2)$. The law of the three-state process is specified by the transition intensities $\alpha_{01}, \alpha_{02}, \alpha_{12}$. There is a way of expressing the intensities $\lambda_1$ and $\lambda_2$ of $N_1$ and $N_2$ such that the laws of $X$ and $N$ correspond to the same probability on $(\Omega, \mathcal{F})$. Thus the same statistical model can be described with $X$ or with $N$.

## 4.3 Statistical models and true probability

So-called objectivist approaches to statistical inference assume that there is a true, generally unknown, probability $\boldsymbol{P}^*$. Frequentists as well as objectivist Bayesians adopt this paradigm while subjectivist Bayesians such as De Finetti (1974) reject it. We adopt the objectivist paradigm which is more suited to answer scientific issues. Statistical inference aims to approach $\boldsymbol{P}^*$ or functionals of $\boldsymbol{P}^*$. Model $\Pi$ is well specified if $\boldsymbol{P}^* \in \Pi$, mis-specified otherwise. If it is well specified there is a $\theta_* \in \Theta$ such that $\boldsymbol{P}^{\theta_*} = \boldsymbol{P}^*$. If we consider a probability $\boldsymbol{P}^\theta$ we may measure its divergence with respect to $\boldsymbol{P}^*$ on a given sigma-field $\mathcal{O}$ by $\boldsymbol{I}(\boldsymbol{P}^\theta | \boldsymbol{P}^*; \mathcal{O})$, and we may choose $\theta$ which minimizes this divergence. We assume that there exist a value $\theta_{opt}$ which minimizes $\boldsymbol{I}(\boldsymbol{P}^\theta | \boldsymbol{P}^*; \mathcal{O})$. We call $\boldsymbol{I}(\boldsymbol{P}^{\theta_{opt}} | \boldsymbol{P}^*; \mathcal{O})$ the misspecification risk of model $\Pi$. Of course if the model is well specified $\boldsymbol{I}(\boldsymbol{P}^\theta | \boldsymbol{P}^*; \mathcal{O})$ is minimized at $\theta_*$ and the misspecification risk is null.

# 5 The likelihood

## 5.1 Definition of the likelihood

Conventionnally most statistical models assume that independently identically distributed (iid) random variables, say $X_i, i = 1, \ldots, n$, are observed. However in case of complex observation schemes the observed random variables become complicated; moreover the same statistical model can be described by different random variables. For instance, in Example 2 the observed random variables are the couples $(\tilde{X}_i, \delta_i)$. However we may also describe the observation by $(\delta_i X_i, \delta_i)$, or in terms of counting processes by $(N_u^i, 0 \leq u \leq C)$, where $(N_u^i = 1_{\{X_i \leq u\}})$. These three descriptions are observationally equivalent, in the sense that they correspond to the same sigma-field, say $\mathcal{O}_i = \sigma(\tilde{X}_i, \delta_i) = \sigma(\delta_i X_i, \delta_i) = \sigma(N_u^i, 0 \leq u \leq C)$. We shall adopt the description of observations in terms of sigma-fields because it is more intrinsic. We shall work with a measure space $(\Omega, \mathcal{F})$ containing all events of interest. For instance the observation of subject $i$, $\mathcal{O}_i$, belongs to $\mathcal{F}$. Saying that observations are iid means that the $\mathcal{O}_i$ are independent, that there is a one-to-one correspondence between $\mathcal{O}_i$ and $\mathcal{O}_{i'}$ and that the restrictions of $\boldsymbol{P}^*$ to $\mathcal{O}_i$, $\boldsymbol{P}^*_{\mathcal{O}_i}$, are the same. We call $\underline{\mathrm{O}}_n$ the global observation: $\underline{\mathrm{O}}_n = \vee_{i=1}^n \mathcal{O}_i$. Since we do not know $\boldsymbol{P}^*$ we may in the first place reduce the search by restricting to a statistical model $\boldsymbol{\Pi}$ and find a $\boldsymbol{P}^\theta \in \boldsymbol{\Pi}$ close to $\boldsymbol{P}^*$, that is, one which minimizes $\boldsymbol{I}(\boldsymbol{P}^\theta | \boldsymbol{P}^*; \mathcal{O}_i)$. We have already given a name to it, $\boldsymbol{P}^{\theta_{opt}}$ but we cannot compute it directly because we do not know $\boldsymbol{P}^*$. The problem is that $\boldsymbol{I}(\boldsymbol{P}^\theta | \boldsymbol{P}^*; \mathcal{O}_i)$ doubly depends on the unknown $\boldsymbol{P}^*$: (i) through the Radon-Nikodym derivative; (ii) through the expectation. Prob-lem (i) can be eliminated by noting that $L_{\mathcal{O}_i}^{\boldsymbol{P}^*/\boldsymbol{P}^\theta} = L_{\mathcal{O}_i}^{\boldsymbol{P}^*/\boldsymbol{P}^0} + L_{\mathcal{O}_i}^{\boldsymbol{P}^0/\boldsymbol{P}^\theta}$. Thus, by taking expectation under $\boldsymbol{P}^*$:

$$\boldsymbol{I}(\boldsymbol{P}^\theta|\boldsymbol{P}^*;\mathcal{O}_i) = \boldsymbol{I}(\boldsymbol{P}^0|\boldsymbol{P}^*;\mathcal{O}_i) - \mathrm{E}_{\boldsymbol{P}^*}(L_{\mathcal{O}_i}^{\boldsymbol{P}^\theta/\boldsymbol{P}^0})$$

Minimizing $\boldsymbol{I}(\boldsymbol{P}^\theta|\boldsymbol{P}^*;\mathcal{O}_i)$ is equivalent to maximizing $\mathrm{E}_{\boldsymbol{P}^*}(L_{\mathcal{O}_i}^{\boldsymbol{P}^\theta/\boldsymbol{P}^0})$. We cannot compute $\mathrm{E}_{\boldsymbol{P}^*}(L_{\mathcal{O}_i}^{\boldsymbol{P}^\theta/\boldsymbol{P}^0})$ but we can estimate it. The law of large numbers tells us that:

$$n^{-1}\sum_{i=1}^n L_{\mathcal{O}_i}^{\boldsymbol{P}^\theta/\boldsymbol{P}^0} \to \mathrm{E}_{\boldsymbol{P}^*}(L_{\mathcal{O}_i}^{\boldsymbol{P}^\theta/\boldsymbol{P}^0}).$$

Thus we may maximize the estimator on the left hand or equivalently the likelihood function $\mathcal{L}_{\underline{\mathrm{O}}_n}^{\boldsymbol{P}^\theta/\boldsymbol{P}^0} = \frac{d\boldsymbol{P}^\theta}{d\boldsymbol{P}^0}|_{\underline{\mathrm{O}}_n}$. The likelihood function is the function $\theta \to \mathcal{L}_{\underline{\mathrm{O}}_n}^{\boldsymbol{P}^\theta/\boldsymbol{P}^0}$. In conclusion the maximum likelihood estimator (MLE) can be considered as an estimator which minimizes a natural estimator of the Kullback-Leiber risk.

## 5.2 Computation of the likelihood

Computation of the likelihood is simple in terms of the probability on the observed $\sigma$-field. The conventional way of specifying a model is in terms of a random variable and a family of distributions $(X, (f_X^\theta(.))_{\theta\in\Theta})$. Then the likelihood for observation $X$ is simply $f_X^\theta(X)$. When the events of interest are represented by stochastic processes in continuous time, it is also possible to define a density and hence a likelihood function. See Feigin (1976) for diffusion processes and Jacod (1975) for counting processes.

Two situations make the computation of the likelihood more complex. The first is when there is incomplete observation of the events of interest. A rather general approach of this problem is through the concept of coarsening, and to make reasonably simple computation the concept of ignorability of the mechanism leading to incomplete data has been promoted (Heitjan and Rubin, 1991). This has been generalized to the stochastic process framework by Commenges and Gégout-Petit (2005) (which also give some general formulas for likelihood calculus). The second situation occurs when the law is described through a conditional probability and the conditioning events are not observed. This is the framework of random effects models. Although conceptually different these two situations lead to the same problem: the likelihood for subject $i$ can be relatively easily computed for a "complete" observation $\mathcal{G}_i$ and the likelihood for the observation $\mathcal{O}_i \subset \mathcal{G}_i$ is the conditional expectation (which derives from the fundamental formula):

$$\mathcal{L}_{\mathcal{O}_i}^{\boldsymbol{P}^\theta/\boldsymbol{P}^0} = \mathrm{E}_{P^0}\left[\mathcal{L}_{\mathcal{G}_i}^{\boldsymbol{P}^\theta/\boldsymbol{P}^0}|\mathcal{O}_i\right].$$

The conditional expectation is expressed as an integral which must be computed numerically in most cases. The only notable exception is the linear mixed effects model where the integral can be analytically computed. For examples of algorithms for non-linear mixed effects see Delyon, Lavielle and Moulines (1999), Guedj, Thiébaut and Commenges (2007) and for general formulas for the likelihood of interval-censored observations of counting processes, Commenges and Gégout-Petit (2007).

## 5.3 Performance of the MLE in terms of Kullback-Leibler risk

We expect good behaviour of the MLE $\hat{\theta}$ when the law of large numbers can be applied and when the number of parameters is not too large. Some cases of unsatisfactory behaviour of the MLE are reported for instance in Le Cam (1990). The properties of the MLE may not be satisfactory when the number of parameters is too large, and especially when it increases with $n$ such as in an example given by Neymann and Scott (1948).

To assess the performance of the MLE we can use a risk which is an extended version of the Kullback-Leibler risk:

$$\mathrm{EKL}(P^{\hat{\theta}}, P^*) = \mathrm{E}_{P^*}(L_{\mathcal{O}_i}^{P^*/P^{\hat{\theta}}}).$$

The difference with the classical Kullback-Leibler risk is that here $P^{\hat{\theta}}$ is random: so $\mathrm{EKL}(P^{\hat{\theta}}, P^*)$ is the expectation of the Kullkack-Leibler divergence between $P^{\hat{\theta}}$ and $P^*$. In parametric models (that is, $\Theta$ is a subset of $\Re^p$) it can be shown (Linhart and Zucchini, 1986; Commenges et al., 2008) that

$$\mathrm{EKL}(P^{\hat{\theta}}, P^*) = \mathrm{E}_{\boldsymbol{P}^*}[L_{\mathcal{X}}^{\boldsymbol{P}^*/\boldsymbol{P}^{\theta_{opt}}}] + \frac{1}{2}n^{-1}\mathrm{Tr}(I^{-1}J) + o(n^{-1}), \qquad (2)$$

where $I$ is the information matrix and $J$ is the variance of the score, both computed in $\theta_{opt}$. This can be nicely interpreted by saying that the risk $\mathrm{EKL}(P^{\hat{\theta}}, P^*)$ is the sum of the misspecification risk $\mathrm{E}_{\boldsymbol{P}^*}[L_{\mathcal{X}}^{\boldsymbol{P}^*/\boldsymbol{P}^{\theta_{opt}}}]$ and the statistical risk $\frac{1}{2}n^{-1}\mathrm{Tr}(I^{-1}J)$. Note in passing that if $\Pi$ is well specified we have $\mathrm{E}_{\boldsymbol{P}^*}[L_{\mathcal{X}}^{\boldsymbol{P}^*/\boldsymbol{P}^{\theta_{opt}}}] = 0$ and $I = J$, and thus $\mathrm{EKL}(P^{\hat{\theta}}, P^*) = \frac{p}{2n} + o(n^{-1})$.

# 6  The penalized likelihood

There is a large literature on the topic: Good and Gaskin (1971), Wahba (1983), O Sullivan (1988), Hastie and Tishirani (1990), Joly and Commenges (1999), Gu and Kim (2002) among others. Penalized likelihood is useful when the statistical model is too large to obtain good estimators, while conventional parametric models appear too rigid. A simple form of the penalized log-likelihood is

$$pl_\kappa(\theta) = \log \mathcal{L}_\mathcal{O}^\theta - \kappa J(\theta).$$

where $J(\theta)$ is a measure of our dislike of $\theta$ and $\kappa$ weights the influence of this measure on the objective function. A classical example is when $\theta = (\alpha(.), \beta)$, where $\alpha(.)$ is a function and $\beta$ is a real parameter. $J(\theta)$ can be chosen as

$$J(\theta) = \int_0^\infty \alpha''(u)^2 du.$$

In this case $J(\theta)$ measures the irregularity of the function $\alpha(.)$. The maximum penalized likelihood estimator (MpLE) $\theta_\kappa^{pl}$ is the value of $\theta$ which maximizes $pl_\kappa(\theta)$. $\kappa$ is often called a smoothing coefficient in the cases where $J(\theta)$ is a measure of the irregularity of a function. More generally, we will call it a meta-parameter. We may generalize the penalized log-likelihood by replacing $\kappa J(\theta)$ by $J(\theta, \kappa)$, where $\kappa$ could be multidimensional. When $\kappa$ varies, this defines a family of estimators, $(\theta_\kappa^{pl}; \kappa \geq 0)$. $\kappa$ may be chosen by cross-validation (see section 8).

There is another way of dealing with the problem of possibly too large statistical models, the so-called sieve estimators. Consider a family of models $(\mathcal{P}_\nu)_{\nu \geq 0}$ where:

$$\mathcal{P}_\nu = (P^\theta; \theta \in \Theta : J(\theta) \leq \nu).$$

For fixed $\nu$, the MLE solves the constrained maximization problem:

$$\max L_{\mathcal{O}}^{\theta}; \text{ subject to } J(\theta) \leq \nu \tag{3}$$

Let us denote $\hat{\theta}_{\nu}$ the MLE. When $\nu$ varies this defines a family of sieve estimators: $(\hat{\theta}_{\nu}; \nu \geq 0)$. $\hat{\theta}_{\nu}$ maximizes the Lagrangian $L_{\mathcal{O}}^{\theta} - \lambda[J(\theta) - \nu]$ for some value of $\lambda$. The Lagrangian superficially looks like the penalized log-likelihood function but an important difference is that here the Lagrange multiplier $\lambda$ is not fixed and is a part of the solution. If the problem is convex the Karush-Kuhn-Tucker conditions are necessary and sufficient. Here these conditions are

$$J(\theta) \leq \nu; \lambda \geq 0; \frac{\partial L_{\mathcal{O}}^{\theta}}{\partial \theta} - \lambda \frac{\partial J(\theta)}{\partial \theta} = 0. \tag{4}$$

It is clear that when the observation $\mathcal{O}$ is fixed, the function $\kappa \to J(\theta_{\kappa}^{pl})$ is a monotone decreasing function. Consider the case where this function is continuous and unbounded (when $\kappa \to 0$). Then for each fixed $\nu$ there exists a value, say $\kappa_{\nu}$, such that $J(\theta_{\kappa_{\nu}}^{pl}) = \nu$. Note that this value depends on $\mathcal{O}$. Now, it is easy to see that $\theta_{\kappa_{\nu}}^{pl}$ satisfies the Karush-Kuhn-Tucker conditions (4), with $\lambda = \kappa_{\nu}$. Thus if we can find the correct $\kappa_{\nu}$ we can solve the constrained maximization problem by maximizing the corresponding penalized likelihood. However, the search for $\kappa_{\nu}$ is not simple and we must remember that the relationship between $\nu$ and $\kappa_{\nu}$ depends on $\mathcal{O}$. A simpler result, deriving from the previous considerations, is:

**Lemma 1 (Penalized and sieves estimators)** *The families $(P^{\theta_{\kappa}^{pl}}; \kappa \geq 0)$ and $(P^{\hat{\theta}_{\nu}}; \nu \geq 0)$ are identical families of estimators.*

The consequence is that since it is easier to solve the unconstrained maximization problem involved in the penalized likelihood approach one should apply this approach in applications. On the other hand it may be easier to develop asymptotic results for sieve estimators (because $\hat{\theta}_\nu$ is a MLE) than for penalized likelihood estimators. One should be able to derive properties of penalized likelihood estimators from those of sieve estimators.

# 7 The hierarchical likelihood

An important class of models arises when we define a potentially observable variable $Y_i$ for each subject and its distribution is given conditional on unobserved quantities. Specifically, let us consider the following model: conditionally on $b^i$, $Y_i$ has a density $f_{Y|b}(.; \theta, b^i)$, where $\theta$ is a vector of parameters of dimension $m$ and $b^i$ are random effects (or parameters) of dimension $K$. The $(Y_i, b^i)$ are i.i.d. Typically $Y_i$ is multivariate of dimension $n_i$. We assume that the $b^i$ have density $f_b(.; \tau)$, where $\tau$ is a parameter. Typically $Y_i$ is observed while $b^i$ is not. This can be made more general for including the case of censored observation of $Y_i$. This is the classical framework of random effects models.

The conventional approach for estimating $\theta$ is to compute the maximum likelihood estimators. Empirical Bayes estimators of the $b^i$ can be computed in a second stage. The likelihood is computed by taking the expectation of the conditional likelihood given the random effect.

$$\mathcal{L}_{\mathcal{O}_i} = \mathrm{E}(\mathcal{L}_{\mathcal{O}_i|b_i}|\mathcal{O}_i).$$

Practically the computation of this conditional expectation involves the integral $\int f^\theta_{Y|b}(Y_i|b) f_b(b) db$. However, the computation of these multiple integrals of dimension $K$ is a daunting task if $K$ is larger than 2 or 3, especially if the likelihood given the random effects is not itself very easy to compute: this is the curse of dimensionality.

For hierarchical generalized linear models the hierarchical likelihood, or h-likelihood, was proposed by Lee and Nelder (1996, 2001); see also Lee, Nelder and Pawitan (2006). The h-likelihood is the joint likelihood of the observations and the (unobserved) random effects. Estimators (here denoted MHLE) of both $\theta$ and $b$ are obtained by maximizing the h-likelihood. In practice this is done by maximizing the h-loglikelihood:

$$hl(\gamma) = L^{\gamma}_{\underline{O}_n} - \sum_{i=1}^{n} \log f_b(b^i; \tau).$$

where $L^{\gamma}_{\underline{O}_n}$ is the loglikelihood for the observation conditional on $b$, and $\gamma = (\theta, b)$ is the set of all the "parameters". $L^{\gamma}_{\underline{O}_n}$ is the likelihood computed as if $b$ were ordinary parameters; in a conventional random-effect approach this would be considered as a conditional likelihood. Often the loglikelihood can be written $L^{\gamma}_{\underline{O}_n} = \sum_i^n \log f(Y_i; \theta, b^i)$. However this formulation is not completely general, because there are interesting cases where observations of the $Y_i$ are censored. So we prefer writing the loglikelihood as $L^{\gamma}_{\underline{O}_n}$ where $\underline{O}_n$ represents the observed $\sigma$-field for $n$ subjects. We note $\hat{\gamma} = (\hat{\theta}, \hat{b})$ the maximum h-likelihood estimators of the parameters for given $\tau$; the latter (meta) parameter can be estimated by profile likelihood. The main interest of this approach is that there is no need to compute multiple integrals. This problem is replaced by that of maximizing $hl(\gamma)$ over $\gamma$, that is, a large number of parameters: this number is equal to $m + nK$. This may be large but special algorithms can be used for generalized linear models.

Therneau and Grambsch (2000) used the same approach for fitting frailty models, calling it a penalized likelihood. It may superficially look like the penalized quasi likelihood of Breslow and Clayton (1993) but this is not the same thing. There is a link with the more conventional penalized likelihood for estimating smooth functions discussed in section 6. The h-likelihood can be considered as a penalized likelihood but with two important differences relative to the conventional one: (i) the problem is parametric; (ii) the number of parameters grows with $n$. Commenges et al. (2008) have proved that the maximum h-likelihood estimators for the fixed parameters are M-estimators (see van der Vaart, 1998). Thus under some regularity conditions they have an asymptotic normal distribution. However, this asymptotic distribution is not in general centered on the true parameter values so that the estimators are biased. In practice the bias can be negligible so that this approach can be interesting in some situations due to its relative numerical simplicity.

## 8 The likelihood cross-validation criterion

An important issue is the choice between different estimators. Two typical situations are : (i) choice of MLE's in different models; (ii) choice of MpLE's with different penalties. If we consider two models $\Pi$ and $\Pi'$ we get two estimators $P^{\hat{\theta}}$ and $P^{\hat{\gamma}}$ of the probability $P^*$ and we may wish to assess which is better: this is the "model choice" issue. A penalized likelihood function produces a family of estimators $(P^{\theta_\kappa^{pl}}; \kappa \geq 0)$ and we may wish to choose the best. Here what we call "the best" estimator is the estimator which minimizes some risk fucntion; in both cases we can use the extended version of the Kullback-Leibler risk already used in section 5:

$$\text{EKL}(P^{\hat{\theta}}, P^*) = \text{E}_{P^*}(L_{\mathcal{O}_i}^{P^*/P^{\hat{\theta}}}).$$

Since $P^*$ is unknown we can first work with $\mathrm{EKL}(P^{\hat{\theta}}, P^0) = \mathrm{E}_{P^*}(L_{\mathcal{O}_i}^{P^0/P^{\hat{\theta}}})$, which is equal, up to a constant, to $\mathrm{EKL}(P^{\hat{\theta}}, P^*)$. Second we can, as usual, replace the expectation under $P^*$ by expectation under the empirical distribution. In parametric models Akaike (1973) has shown that an estimator of $\mathrm{EKL}(P^{\hat{\theta}}, P^0)$ was $-n^{-1}(\mathcal{L}_{\mathcal{O}_i}^{P^{\hat{\theta}}/P^0} - p)$. The AIC criterion can be deduced by multiplying by $2n$. The result can be used to estimate the difference of risks between two estimators in parametric models $\Delta(P^{\hat{\theta}}, P^{\hat{\gamma}}) = \mathrm{EKL}(P^{\hat{\theta}}, P^*) - \mathrm{EKL}(P^{\hat{\gamma}}, P^*)$ by the statistic $D(P^{\hat{\theta}}, P^{\hat{\gamma}}) = (1/2n)(AIC(P^{\hat{\theta}}) - AIC(P^{\hat{\gamma}}))$ and a more refined analysis of the difference of risks can be developed, as in Commenges et al. (2008).

## 9 Link withe MAP estimator

One important issue is the relationship between the three likelihoods considered here and the Bayesian approach. The question arises because it seems that these three likelihoods can be identified with the numerator of *a posteriori* distributions with particular priors. Thus MLE, MpLE and MHLE could be identified with the maximum a posteriori (MAP) estimators with the corresponding priors. However, this relationship depends on the parametrization. Thus the MLE is identical to the MAP using a flat prior for the parameters; if we change the parametrization, the flat prior on the new parameters does not correspond to the flat prior on the original parameters, as was already noticed by Fisher (1922). This apparent paradox led Jeffreys to propose a prior invariant for parametrization (Jeffreys, 1961), known as Jeffrey's prior. However the MAP with Jeffreys's prior is no longer identical to the MLE when Jeffreys's prior is not flat. For instance for the parameter of binomial trial Jeffreys's prior is $1/\sqrt{p(1-p)}$. Adding the logarithm of this term to the loglikelihood shifts the maximum away from 0.5. Moreover it is questionable whether this invariance property can be identified with a non-informativeness character of this prior (for a review on the choice of priors, see Kass and Wasserman, 1996).

In the Bayesian paradigm, rather than considering estimators based on maximization of some expression such as the likelihood or posterior density, it is common to attempt to summarize the statistical inferences by using qunatiles of the posterior distribution, such as the median, or expectations with respect to the posterior. While such expectation may be more satisfactory, they typically involve multiple integrals which are hard to compute: computations are mostly being done with the MCMC algorithm. Maximization methods have the advantage of being potentially easier in the case where multiple integrals can be avoided. There are also approximate Bayesian methods which yield the a posteriori marginal distribution by approximating some of the multiple integrals by Laplace approximation, which in turn involves a maximization problem: Rue et al., (2008) claim that this approach is much faster than the MCMC algorithm.

## Conclusion

The Kolmogorov representation of a statistical experiment has to be taken seriously if we want to have a deep understanding of what a statistical model is. The Kulback-Leibler risk is underlying most of the reflexions about likeli-hood, as was clearly seen by Akaike (1973). Finally the link with the Bayesian approach should be explored more thoroughly than could done in this paper. The MLE and MAP estimators are the same if, in a given paramtrization, the prior used for the MAP is flat. However, this identity does not resist to a reparametrization. Similar remarks hold for the link between penalized likelihood and MAP.

**Acknowledgement**. I would like to thank Anne Gégout-Petit for helpful comments on the manuscript.